\documentclass[12pt]{amsart}
\usepackage{amsfonts, amssymb,amsmath, amsthm, amsxtra, latexsym, amscd}
\usepackage[latin1]{inputenc} 

\newtheorem{theo+}    {Theorem}      [section]
\newtheorem{prop+}  [theo+]  {Proposition}
\newtheorem{coro+}  [theo+]  {Corollary}
\newtheorem{lemm+}  [theo+]  {Lemma}
\newtheorem{deep+}  [theo+]  {Deep Result}
\newtheorem{fact+}  [theo+]  {Fact}
\theoremstyle{definition}
\newtheorem{exam+}  [theo+]  {Example}
\newtheorem{rema+}  [theo+]  {Remark}
\newtheorem{defi+}  [theo+]  {Definition}
\newtheorem{xca+}[theo+]{Exercise}

 \numberwithin{equation}{section}

\def\beqn{\begin{equation}}
\def\eeqn{\end{equation}}

\def\epf{\qed \enddemo}

\def\a{\alpha}

\def\epsi{\epsilon}
\def\Claminv2{|C(\Lambda)|^{-2}}

\def\varepsi{\varepsilon}
\def\lam{\lambda}

\def\blam{\underline{\bold \lambda}}

\def\Ome{\Omega}

\def\Aa2D{A^{\a,2}(D)}
\def\bAa2D{\overline{A^{\a,2}(D)}}
\def\Ab2D{A^{\beta,2}(D)}
\def\bAb2D{\overline{A^{\beta,2}(D)}}

\def\Norm#1_#2{\Vert#1\Vert_{#2}}

\def\2pd#1#2{\frac{\partial^2 #1}{\partial #2^2}}
\def\p11d#1#2#3{\frac{\partial^2 #1}{  \partial #2\partial #3  }}

\def\Claminv2{|C(\Lambda)|^{-2}}

\def\varepsi{\varepsilon}
\def\sig{\sigma}

\def\lam{\lambda}

\def\tanh{\operatorname{tanh}}

\def\tanh{\operatorname{tanh}}

\def\Aa2D{A^{\a,2}(D)}
\def\bAa2D{\overline{A^{\a,2}(D)}}
\def\Ab2D{A^{\beta,2}(D)}
\def\bAb2D{\overline{A^{\beta,2}(D)}}

\def\m{\underline{\mathbf m}}
\def\n{\underline{\mathbf n}}

\def\ub1#1{\underline{\mathbf 1^{#1}}}

\def\h-g-o-p{hypergeometric orthogonal polynomial }
\def\h-g-o-ps{hypergeometric orthogonal polynomials }

\def\bc{\mathbb C}

\def\nat0{\mathbb Z_{\ge 0}} 

\def\bpf{\begin{proof}}
\def\epf{\end{proof}}
\def\beq{\begin{equation}}
\def\eeq{\end{equation}}

\def\bc{\mathbb C}

\def\draft{\centerline{(Draft {\the \day}/{\the\month} \the \year.)}}

\begin{document}

\def\phila{\phi_{\blam}}
\def\vepm{\varepsilon_{\m, \nu}(\blam)}
\def\vepmp{\varepsilon_{\m^\prime, \nu}(\blam)}
\def\kam{\kappa_{\m, \nu}(\blam)}
\def\nutnu{\pi({\nu})\otimes \overline{\pi(\nu)}}
\def\pia2ta2{\pi({\frac a2})\otimes \overline{\pi(\frac a2)}}
\def\Hnu{\mathcal H_{\nu}}
\def\Fnu{\mathcal F_{\nu}}
\def\Pm{\mathcal P_{\m}}
\def\PV{\mathcal P(V_{\bc})}
\def\Pa{\mathcal P(\mathfrak a)}
\def\Pn{\mathcal P_{\n}}
\def\cL#1nu{\mathcal L_{{#1}, \nu}}
\def\cL#1a2{\mathcal L_{{#1}, \frac a2}}
\def\E#1nu{E_{{#1}, \nu}}
\def\E#1a2{E_{{#1}, \frac a2}}

\def\Dc{D_{\mathbb C}}
\def\Vc{V_{\mathbb C}}
\def\Bnu{B_{\nu}} 
\def\el{e_{\blam}} 
\def\bnul{b_{\nu}(\lam)} 
\def\GGa{\Gamma_{\Ome}} 
\def\fc#1#2{\frac{#1}{#2}} 
\def\SS{\mathcal S}

\title[Non-symmetric Jacobi and Wilson type polynomials]
{Non-symmetric Jacobi and Wilson type
polynomials}

\author{Lizhong Peng}
\address{LMAM, School of Mathematical Sciences, Peking University, Peking 100871. China}
\email{lzpeng@math.pku.edu.cn}

\author{Genkai Zhang}
\address{Department of Mathematics, Chalmers University of Technology
and G\"o{}teborg University, S-412 96 G\"o{}teborg, Sweden}
\email{genkai@math.chalmers.se}

\thanks{Research supported by Swedish Research Council
(VR), NNSF of China No. 10471002, and the Swedish Foundation for
international cooperation in higher education (STINT) and
Mathematical Center of Chinese Ministry of Education. }

\begin{abstract}
Consider  a root system  of type $BC_1$ on
the real line $\mathbb R$ with general positive multiplicities.
The Cherednik-Opdam transform defines
a unitary operator from an $L^2$-space on $\mathbb R$
to a $L^2$-space of $\mathbb C^2$-valued functions
on $\mathbb R^+$ with the Harish-Chandra measure $|c(\lam)|^{-2}d\lam$.
 By introducing
a weight function of the form $\cosh^{-\sig}(t)\tanh^{2k} t$ 
on $\mathbb R$ we
find an orthogonal basis 
for the $L^2$-space  on $\mathbb R$ consisting of
even and odd functions expressed in terms 
of the  Jacobi polynomials (for each fixed $\sig$ and $k$). We find a
Rodrigues type formula for the  functions in terms of the
Cherednik operator. We compute explicitly their Cherednik-Opdam
transforms.
We  discover thus a new family of $\mathbb C^2$-valued orthogonal
polynomials.
In the special case when $k=0$ the even polynomials become Wilson polynomials, and the
corresponding result  was proved earlier by Koornwinder.
\end{abstract}


\maketitle \baselineskip 1.40pc

\section{Introduction and Main result}

The Hermite polynomials and functions on the Euclidean space $\mathbb R$
play an important role in Fourier analysis and
in the theory of special functions. The Hermite functions
are a  product of the Gaussian function and the Hermite polynomials,
 they form an orthogonal basis of the space $L^2(\mathbb R)$ and their
Fourier transforms are also Hermite type functions.
A natural  generalization of the Fourier transform is the Jacobi transform
for even functions on $\mathbb R$.
In his paper
\cite{Koornwinder-lnm} Koornwinder
computed the Jacobi transforms of certain even orthogonal function of Jacobi type
and proved
that they are,  up to some factor of  Gamma functions,
 the Wilson hypergeometric orthogonal polynomials.
The Jacobi transform is, for special parameters,
the Harish-Chandra spherical transform in rank  one.
Recently \cite{Opdam-acta} Opdam introduced
a generalization of the Harish-Chandra spherical
transform for general non-symmetric functions on
any root system $R$
in
 $\mathbb R^r$
with general positive root multiplicities,
  and proved
the corresponding  transform is unitary from
 certain $L^2$-space
to a space of vector-valued functions
on the positive Weyl chamber in
$\mathbb R^r$ with the Harish-Chandra measure
$|c(\lam)|^{-2}d\lam$.
The transform is
 also called the Cherednik-Opdam transform.
When restricted
to symmetric (i.e. Weyl group invariant) functions on
 root systems of a non-compact symmetric
space $G/K$, the transform reduces to the  spherical
transform
for $K$-invariant functions.

There appears thus a natural  question,
namely
 to find an orthogonal basis on
$\mathbb R^r$ and to compute its Cherednik-Opdam transform.
Motivated
by the study of Berezin transform
\cite{gkz-bere-rbsd},  \cite{Dijk-pevzner},
and branching rule \cite{gz-br2},
Zhang introduces \cite{gz-sphcan}
certain weight function, also called canonical weight $w_\sig(t)$ on root
system of Type BC of arbitrary rank, computes the spherical
transform $\widetilde w_\sig(\lambda)$ and proves that the
(Weyl group-invariant) Jacobi polynomials $P^J$
multiplied by the weight function, have their spherical transform
being multi-variable  Wilson polynomials $P^W$  \cite{van-Dijen-tams} multiplied
by the spherical transform $\widetilde w_\sig(\lambda)$
of $w_\sig$, thus
generalizing the classical theory of 
Fourier transform of Gaussian and Hermite functions. The rank one case
has been done earlier by Koornwinder \cite{Koornwinder-lnm}
where the spherical transform is the classical
Jacobi transform.
In the present paper we will study
the non-symmetric analogue of the above result
of Koornwinder in rank one in the setup
of Opdam-Cherednik transform.

Let $R=\{\pm 2\varepsi, \pm 4\varepsi\}\subset \mathbb R \varepsi$ be a root system on $\mathbb R$ of type BC with
general positive multiplicities.
Consider the $L^2$-space $L^2(\mathbb R, d\mu)$
with $d\mu$ given in   (\ref{eq:dmu}). Consider the weight function
 $w_{\sig, k}(t)=(\cosh^{-\sig} t)\tanh^{2k}t$
for fixed $\sig$ and  integer $k$
and the system $\{
w_{\sig, k}(t) (\tanh t)^n\}_{n=0}^\infty$; the latter
forms a basis of the $L^2$-space and we find an orthogonal basis of
the form $w_{\sig, k}(t) p_n(\tanh t)
$ where $p_n$ can be written 
in terms of Jacobi polynomial of degree $n$ (see Lemma 3.1).
 We find a Rodrigues type formula for the functions
in Theorem 3.5, and their Cherednik-Opdam transforms are computed
in Theorem 4.1 and Theorem 4.3, written in terms
of certain ${}_3F_2$-hypergeometric series. The
functions and their transform can
be viewed as non-symmetric analogues
of the Jacobi and Wilson polynomials explaining
 the title of this paper.

Part of the work was done when one of the authors, Genkai Zhang,
was visiting the School of Mathematical Sciences at Peking University,
he would like to thank the school for its
hospitality.

\section{Preliminaries. Cherednik-Opdam transform}

We recall briefly the Plancherel formula for non-symmetric functions
proved by Opdam \cite{Opdam-acta} applied to the
root system of type $BC_1$ (in which case it can also be deduced
from the known Plancherel formula for the Jacobi transform
for symmetric functions).
We use the same notations as in \cite{gz-sphcan}.
Let $R=\{\pm 2\varepsi,
\pm 4\varepsi
\}$ be a root system of type $BC_1$ with multiplicities
$k_{2\varepsi}=b$ and $k_{4\varepsi}=\frac{\iota}2$ respectively.
 We normalize $\varepsi$
as a unit vector and identify $\mathbb R \varepsi$
with $\mathbb R$ and with its dual space.
 The half sum of positive
roots is then $\rho= b+\iota$.

Let $d\mu$ be the measure
\begin{equation}
  \label{eq:dmu}
d\mu(t)=\prod_{\a=2\varepsi, 4\epsi}|2\sinh (\frac
12\a(t))|^{2k_{\a}}dt
=2^{2b+\iota}|\sinh t|^{2b} |\sinh(2t)|^{\iota}dt
\end{equation}
on $\mathbb R$, and let  $L^2(\mathbb R, d\mu)$ be the corresponding
 $L^2$-space. The Weyl group  $\{\pm 1\}$ acts on $L^2(\mathbb R, d\mu)$
 via $f(t)\mapsto f(\pm t)$, and under this action the space
is decomposed as
$$
L^2(\mathbb R, d\mu)
=L^2(\mathbb R, d\mu)_{1}\oplus
L^2(\mathbb R, d\mu)_{-1},
$$
a sum of subspaces of even and respectively odd functions.

Let
\begin{equation}
D:=\partial  +2\iota\frac{1}{1-e^{-4t}}(1-s)
+2b\frac{1}{1-e^{-2t}}(1-s) -\rho,
\end{equation}
be the Cherednik operator on  $L^2(\mathbb R, d\mu)$.
Here $s f(t)=f(-t)$ is the reflection.
The eigenvalue problem of $D$ is solved in \cite{Opdam-acta}; the function
\begin{equation*}
\begin{split}
G(\lambda, t) =&{}_2F_1(\frac{\lambda + \rho}2, \frac{-\lambda +
\rho}2; \frac{ 1+\iota}2 +b, -\sinh^2 t)\\
+& \frac 1{\lambda + \rho} \sinh(2t) {}_2F_1^\prime(\frac{\lambda
+ \rho}2, \frac{-\lambda + \rho}2; \frac{ 1+\iota}2 +b, -\sinh^2
t)
\end{split}
\end{equation*}
is an eigenfunction of $D$,
\begin{equation}\label{eig-eq}
DG(\lambda, t)=\lam G(\lambda, t).
\end{equation}
Here ${}_2F_1(a_1, a_2; b, x)$
 is the Gauss hypergeometric function.
For later purpose we recall that
the  hypergeometric function is defined by
$$
{}_pF_q(a_1, \cdots,  a_p; b_1, \cdots, b_q, x)=\sum_{m=0}^\infty \frac{(a_1)_m \cdots (a_p)_m}{(b_1)_m
\cdots (b_q)_m}\frac{x^m}{m!}
$$
where $(a)_m=a(a+1)\cdots (a+m-1)$ is
the Pochhammer symbol;
see e. g. \cite{Erdelyi-1}.

The Cherednik-Opdam transform is defined by, for $f\in C_0^\infty (\mathbb R)$,
as a $\mathbb C^2$-valued function,
$\mathcal Ff(\lambda )=(\mathcal F_{ 1}f(\lambda), \mathcal F_{-1}f(\lambda))$ in $\lambda$,
with
$$
\mathcal F_{\pm 1}f(\lambda)
=\int_{\mathbb R} f(t)G(\lambda, \pm t)d\mu(t).
$$
Let $d\hat \mu(\lam)$ be the measure
$$
d\hat \mu(\lam)=(2\pi)^{-1}\frac{c_{-1}(\rho)^2}{c(\lam)c(-\lam)}d\lam
$$
on the imaginary half axis $i\mathbb R^+$, with
$$
c(\lam)=\frac{ \Gamma(\lam ) \Gamma(\frac \lam 2 +b) }
{\Gamma(\lam +b ) \Gamma(\frac \lam 2 +b +\frac{\iota}2 )
}
$$
and
$$
c_{-1}(\lam)=\frac{ \Gamma(\lam +1 ) \Gamma(\frac \lam 2 +b +1) }
{\Gamma(\lam +b +1 ) \Gamma(\frac \lam 2 +b +\frac{\iota}2 +1).
}
$$
Then the transform $\mathcal F$ extends to a unitary
operator from $L^2(\mathbb R, d\mu(t))$
onto

\noindent $L^2(i\mathbb R^+, d\hat \mu(\lam))\otimes \mathbb C^2$,
namely
$$
\Vert f\Vert^2 =\Vert \mathcal Ff \Vert^2_{
L^2(i\mathbb R^+, d\hat \mu(\lam) )\otimes \mathbb C^2 }.
$$
When $f$ is an even function
$\mathcal F_{1}f(\lam)=\mathcal F_{-1}f(\lam):=\widetilde f(\lam)$
is the spherical transform of $f$;
see \cite{Opdam-acta} and
\cite{Cherednik-imrn97}.

\section{Jacobi type functions
and Rodrigues type formulas}

In this section
we will find an orthogonal basis
of the space $L^2(\mathbb R, d\mu)$ consisting
of Jacobi type functions
and find certain  Rodrigues type formulas.
 We denote
\begin{equation}
  \label{eq:weigt}
w_{\sig, k}(t):=(\cosh^{-\sig}t) \tanh^{2k}t,
\qquad w_{\sig}(t):=w_{\sig, 0}(t)=\cosh^{-\sig}t.
\end{equation}
The function $w_\sig$  will play the role of the Gaussian function.

Let $P_n^{\alpha, \beta}$ be the Jacobi polynomials \cite{Erdelyi-1}.
We define the even and odd Jacobi-type functions as
follows,
\begin{equation}
Q_{n, 1}^{(k)}(t)=w_{\sig, k}(t) P_n^{\sig_0, \delta_0+2k}
(2\tanh^2
t -1), \qquad n=0, 1, \cdots, \end{equation}
and respectively
\begin{equation}
Q_{n, -1}^{(k)}(t)= w_{\sig, k}(t) P_n^{\sig_0, \delta_1+2k}
(2\tanh^2 t -1) \tanh t, \qquad n=0, 1, \cdots
\end{equation}
where
\begin{equation}
\label{sig-delta-01}
\sig_0:=\sig -(\iota +b +1), \quad
\delta_0:=\frac{\iota -1}2+b, \quad
 \delta_1:=  \delta_0 +1=\frac{\iota +1}2+b .
\end{equation}
\begin{lemm+} Suppose $\sig>\iota +b$.
The functions $\{Q_{n, 1}^{(k)}(t)\}$ and  $\{Q_{n,
-1}^{(k)}(t)\}$ form orthogonal bases for $L^2(\mathbb R,
d\mu(t))_1$ and $ L^2(\mathbb R, d\mu(t))_{-1}$ respectively.
Their norms in the $L^2$-space are given by
$$
\Vert Q^{(k)}_{n, 1}\Vert^2=2^{2(\iota +b)}\frac{\Gamma(n+\sig
-\iota-b-1) \Gamma(n+2k+\frac{\iota +2b +1}2)} {n! (2n +2k +\sig
-\frac{\iota +1}2) \Gamma(n +2k+\sig -\frac{\iota +1}2) }$$ and
$$
\Vert Q^{(k)}_{n, -1}\Vert^2= 2^{2(\iota +b)}\frac{\Gamma(n+\sig
-\iota-b-1) \Gamma(n+2k+\frac{\iota +2b +1}2 +1)} {n! (2n +2k+\sig
-\frac{\iota -1}2) \Gamma(n +2k+\sig -\frac{\iota -1}2) }.
$$
\end{lemm+}

\begin{proof} We perform first the change of variables
$t\in \mathbb R\mapsto z=\tanh t\in (-1, 1)$. Then
the measure $d\mu(t)$, written in the variable $z$, is
(with some abuse of notation)
$$
d\mu(z):=d\mu(t)=2^{2b+2l}|\sinh^{\iota +2b} t \cosh^{\iota} t|dt
=2^{2b+2l}|z^{\iota +2b} (1-z^2)^{-(\iota +b +1)}|dz.
$$
Thus the map $g(t)\mapsto f(z):=f(\tanh t)=g(t)$  is  unitary
from the space $L^2(\mathbb R, d\mu(t))_{\pm 1}$
onto $L^2((-1, 1), d\mu(z))_{\pm 1}$ and we will thus identify
the two spaces.
We consider further the space $L^2((-1,1), \mu_{\sig, k})$
on $(-1, 1)$ with the measure 
$$
d\mu_{\sig, k}(x):=2^{3(\iota +b)-\sig -\frac{\iota +2b-1}2 -2k}
(1-x)^ {\sig -(\iota +b +1)}
(1+x)^{\frac{ \iota +2b-1}2+2k} dx,$$
and the  operator $U: L^2((-1, 1), d\mu(z))_1
=L^2(\mathbb R, d\mu(t))_1 \mapsto
L^2((-1,1), \mu_{\sig, k})$ defined by
\begin{equation*}
\begin{split}
F(x):&=(Uf)(x):= f(z)(1-z^2)^{-\frac{\sig}2} z^{-2k}\\
 &=f(\tanh t)w_{\sig, k}(t)^{-1}, \quad x= 2z^2-1=2\tanh^2 t-1, z=\tanh t;
\end{split}
\end{equation*}
its inverse is given by
$$
f(\tanh t)=f(\tanh(-t))=F(2\tanh^2 t-1)  w_{\sig, k}(t).
$$
Then
\begin{equation*}
  \begin{split}
\Vert f\Vert^2&=\int_{-1}^1 |f(z)|^2 d\mu(z)
=2\int_{0}^1 |f(z)|^2 d\mu(z)\\
&=2^{3(\iota +b)-\sig -\frac{\iota +2b-1}2 -2k}
\int_{-1}^1 |F(x)|^2(1-x)^ {\sig -(\iota +b +1)}
(1+x)^{\frac{ \iota +2b-1}2+2k} dx,
  \end{split}
\end{equation*}
namely the map $f\to F=Uf$ is an isometry, and it is clearly onto. Now the Jacobi polynomials
 form an orthogonal basis of the later space,
 whose norm can be computed
by using the result
that
$$
\int_{-1}^1(1-x)^{\alpha}
(1+x)^{\beta}P_{m}^{(\alpha, \beta)}(x)
P_{n}^{(\alpha, \beta)}(x)dx
=\delta_{m,n}\frac{2^{\alpha +\beta +1}\Gamma(n +\alpha +1)
\Gamma(n +\beta +1)}
{n!(2n +\alpha +\beta +1)\Gamma(n +\alpha+\beta +1)};
$$
see e.g.  \cite{Erdelyi-1}.
Our
result for
$\{Q^{(k)}_{n, 1}\}$ follows immediately.
The claim for  $L^2(\mathbb R, \mu(t)dt)_{-1
}$ follows
by considering the map  $f(t)\mapsto \frac{f(z)}{z}$,  $z=\tanh t$
and by similar computation as above; we  omit the elementary
computations.
\end{proof}

We shall now find certain Rodrigues type formula for the functions
$\{Q^{(k)}_{n, \pm1}\}$ in terms of the weight $w_\sig(t)=w_{\sig, 0}(t)$
and the Cherednik operator $D$. We consider first the case $k=0$
for which we have rather compact formulas,
and the case for general $k\ge 0$ is somewhat different.

We recall the formula for
the  spherical transform of the weight function
$w_{\sig}(t)$ and
certain Bernstein-Sato type formula
proved in \cite{gz-sphcan}, formulated here
slightly differently using the Pochhammer symbol.
\begin{lemm+} \label{bern}
\begin{itemize}
\item[(i)] Let $B_{m, \sig}(x)$ and
$b_{m, \sig}$ be defined by
$$
B_{m, \sig}(x)=\prod_{\epsi =\pm}\left(\frac {\sig-\rho +\epsi x}2\right)_{m},
\quad
b_{m, \sig}=\left(\frac{\sig}2\right)_m\left(\frac{\sig +1-\iota}2\right)_{m}.
$$
Then we have a Bernstein-Sato type formula for
the weight function $w_\sig(t)=w_{\sig}(t)$,
\begin{equation}\label{shift-1}
B_{m, \sig}(D)w_{\sig}(t)
=b_{m, \sig} w_{\sig+2m}(t).
\end{equation}
\item[(ii)] Suppose $\sig >2(\iota +b)$.
The spherical transform $\widetilde{w_{\sig}}=\mathcal F_{\pm}{w_{\sig}}$ of $w_{\sig}(t)$ is
given by
$$
\widetilde{w_{\sig}}(\lam)=2^{2\iota +2b}
\Gamma(\frac{\iota+1 +2b }2) \frac{\Gamma(\frac{\sig}2 +\iota +b))}
{\Gamma (\frac{\sig +1-\iota}2)}\prod_{\epsilon=\pm}
\frac{\Gamma(\frac 12(\sig-\rho)+\epsilon \frac 12 \lam)}{
\Gamma(\frac 12(\sig-\rho)+\epsilon \frac 12(\iota+b ))},
\quad \lam \in i \mathbb R.
$$
\end{itemize}
\end{lemm+}

\begin{rema+}
Note that the conditions on $\sig$ in Lemma 3.1
and in Lemma 3.2(ii) are different.
The weaker condition  $\sig >(\iota +b)$
in Lemma 3.1 is necessary and sufficient
for the functions $w_{\sig}$ in $L^2(\mathbb R, d\mu(t))$
whereas the condition  $\sig >2(\iota +b)$
is so that  the function
$w_{\sig}$ is $L^1(\mathbb R, d\mu(t))$. Lemma 3.2 (ii) and
Theorem
4.1 below can be proved to true for
the weaker condition $\sig >(\iota +b)$ using  some certain arguments
on analytic continuation and we will not present
them there.
\end{rema+}

Applying the operator $D+\rho$ on the above identity
and using the fact that
$$
(D+\rho)\cosh^{-\sig-2m}t
= (-\sig -2m)(\cosh^{-\sig-2m}t) \tanh t$$
 we find
a similar  Bernstein-Sato type formula.
\begin{lemm+} \label{var-bern}
Let $\sig \in \mathbb C$. The following
formula holds
\begin{equation*}
\begin{split}
&(D+\rho)\prod_{\epsi =\pm}\left(\frac {\sig-\rho +\epsi
D}2\right)_{m} w_{\sig}(t) \\
= &(-\sig)\left(\frac{\sig}2 +1\right)_m
\left(\frac{\sig+1-\iota}2\right)_m w_{\sig+2m}(t)\tanh t.
\end{split}
\end{equation*}
\end{lemm+}
This lemma is quite
similar to (\ref{shift-1})
except that there is now an extra factor of $D+\rho$
in the left hand side and a  factor $\tanh t$
in the right hand, and  that the  factor $(\frac \sig 2)_m$
is changed to  $(\frac \sig 2 +1)_m$. This
difference will appear in all the formulas
below between the even polynomials
$Q_{n, 1}^{(k)}$ and 
the odd polynomials $Q_{n, -1}^{(k)}$.

The orthogonal functions $Q^{(0)}_{n, \pm 1}(t)$ can now be obtained
by some polynomials of $D$ acting on the weight function
$w_{\sig}$.

\begin{theo+}
 \label{Rodrigue}
 Let $\mathcal Q^{(0)}_{n, 1}(x)$
and
$\mathcal  Q^{(0)}_{n, -1}(x)$ be the following
polynomial
\begin{equation*}
\begin{split}
&\qquad \mathcal Q^{(0)}_{n, 1}(x)=\frac{(\sig_0 +1)_n}{n!} \times \\
&{}_4F_3(-n, n +\sig_0 +\delta_0 +1, \frac{\sig-\rho +x}2,
\frac{\sig-\rho -x}2; \sig_0+1, \frac{\sig+1-\iota}2,
\frac{\sig}2, 1),
\end{split}
\end{equation*}
\begin{equation*}
\begin{split}
&\qquad \mathcal Q^{(0)}_{n, -1}(x)=\frac{(\sig_0 +1)_n}{n!(-\sig)}
(\rho +x) \times \\
& {}_4F_3 (-n, n +\sig_0 +\delta_1 +1, \frac{\sig-\rho +x}2,
\frac{\sig-\rho -x}2; \sig_0+1, \frac{\sig+1-\iota}2, \frac{\sig}2
+1, 1).
\end{split}
\end{equation*}
We have the following Rodrigue's type formulas
$$
Q^{(0)}_{n, 1}(t)=\mathcal Q^{(0)}_{n, 1}(D) w_{\sig}(t),
$$
and
$$
Q^{(0)}_{n, -1}(t)=\mathcal Q^{(0)}_{n, -1}(D)w_{\sig}(t).
$$
\end{theo+}
\begin{proof}
 The function $Q^{(0)}_{n, 1}$
is
\begin{equation*}
\begin{split}
Q^{(0)}_{n, 1}(t)&= \frac{(\sig_0 +1)_n}{n!} {}_2F_1 (-n, n+\sig_0
+\delta_0 +1; \sig_0 +1;
\cosh^{-2}t)\cosh^{-\sig}t\\
&= \frac{(\sig_0 +1)_n}{n!} \sum_{m=0}^n \frac{(-n)_m (n+\sig_0+\delta_0 +1)_m}
{(\sig_0 +1)_m} \cosh^{-\sig-2m}t.
\end{split}
\end{equation*}
We rewrite Lemma 3.2 as
\begin{equation}
  \label{eq:bern-var2}
w_{\sig+ 2m}(t)=\cosh^{-\sig-2m}t = 
\frac{B_{m, \sig}(D)\cosh^{-\sig}t}{b_{m, \sig}}.
\end{equation}
Our first formula follows then immediately by rewriting
the sum as a hypergeometric function.

Similarly,
\begin{equation*}
\begin{split}
Q^{(0)}_{n, -1}(t)& = \frac{(\sig_0 +1)_n}{n!}{}_2F_1 (-n,
n+\sig_0 +\delta_1 +1; \sig_0 +1;
\cosh^{-2}t) (\cosh^{-\sig}t) \tanh t\\
&= \frac{(\sig_0 +1)_n}{n!} \sum_{m=0}^n \frac{(-n)_m
(n+\sig_0+\delta_1 +1)_m}
{(\sig_0 +1)_m}(\cosh^{-\sig-2m}t) \tanh t,
\end{split}
\end{equation*}
and by Lemma 3.4 each term $(\cosh^{-\sig-2m}t) \tanh t$ in the sum
can be rewritten
as a polynomial of $D$ acting on $w_\sig (t)$. This
proves our claim.
\end{proof}

We consider now the general case of $k\ge 0$.
We shall first generalize Lemma 3.2 and Lemma 3.4 
and find a Bernstein Sato type
formula expressing the  functions
 $w_{\sig+2m, k}(t)$
and  $w_{\sig+2m, k}(t)\tanh t$
as a polynomial   of the operator $D$ acting
on $w_{\sig}(t)=\cosh^{-\sig}t$.

\begin{prop+} \label{bern-k} Let
$L_{k, m}(x)$ and $M_{k, m}(x)$ be the polynomials
\begin{equation}\label{L-km}
L_{k, m}(x)={}_3F_2(-k, \frac{\sig-\rho +}2,\frac{\sig-\rho
-D}2;  \frac\sig 2 +m, \frac{\sig +1-\iota}2, 1)
\end{equation}
\begin{equation}\label{M-km}
M_{k, m}(x)
={}_3F_2(-k, \frac{\sig-\rho +}2,\frac{\sig-\rho
-D}2;  \frac\sig 2 +1 +m, \frac{\sig +1-\iota}2, 1).
\end{equation}
Then the following Bernstein-Sato type formulas hold
\begin{equation}
B_{m, \sig}(D)L_{k, m}(D)
w_{\sig}(t)
=b_{m, \sig} w_{\sig+2m, k}(t),
\end{equation}
\begin{equation}
\quad (D+\rho)
B_{m, \sig}(D) 
M_{k, m}(D)
w_{\sig}(t)=(-\sig)(\frac \sig 2+1)_m
(\frac {\sig +1-\iota}2+1)_m
  w_{\sig+2m, k}(t) \tanh t.
\end{equation}
\end{prop+}

\begin{proof} We write the
weight function
$w_{\sig+2m, k}(t)$ as
\begin{equation*}
\begin{split}
&\quad w_{\sig+2m, k}(t)\cosh^{-(\sig+2m)}\tanh^{2k}t
=\cosh^{-(\sig+2m)}(1-\cosh^{-2})^k\\
&=\sum_{j=0}^k(-1)^j\binom k j \cosh^{-(\sig+2m+2k)}t\\
&=\sum_{j=0}^k\frac{(-k)_j}{j!}
\cosh^{-(\sig+2m+2j)}t.
\end{split}
\end{equation*}
By the Lemma 3.2 each term $\cosh^{-(\sig+2m+2j)}$ is
\begin{equation*}
\begin{split}
&\qquad \cosh^{-(\sig+2m+2j)}t
=\frac{(\frac{\sig -\rho + D}2)_{m+j}
(\frac{\sig -\rho - D}2)_{m+j}}
{(\frac{\sig}2)_{m+j} (\frac{\sig+1-\iota}2)_{m+j}
}\cosh^{-\sig}t\\
&=\frac{(\frac{\sig -\rho + D}2)_{m}
(\frac{\sig -\rho - D}2)_{m}}
{(\frac{\sig}2)_{m} (\frac{\sig+1-\iota}2)_{m}
}
\frac{(\frac{\sig -\rho + D}2 +m)_{j}
(\frac{\sig -\rho - D}2+m)_{j}}
{(\frac{\sig}2+m)_{j} (\frac{\sig+1-\iota}2+m)_{j}}
\cosh^{-\sig}t.
\end{split}
\end{equation*}
Our first result now 
follows by writing the sum as a hypergeometric series. The
second formula is proved by similar computations
using Lemma 3.4.
\end{proof}
\section{ Cherednik-Opdam  transform
of the Jacobi type functions}

We compute now the Cherednik-Opdam transform of the orthogonal
functions

\noindent $\{Q^{k}_{n, \pm 1}\}$. Consider the case $k=0$ first.

\begin{theo+} Suppose $\sig > 2(\iota +b) $.  The Cherednik-Opdam transform of the
Jacobi type functions $Q^{(0)}_{n, 1}$ and $Q^{(0)}_{n, -1}$ are given
by
\begin{equation*}
\begin{split}
&\qquad \mathcal F_{\pm 1}(Q^{(0)}_{n, 1})=\frac{(\sig_0 +1)_n}{n!}
\widetilde{w_{\sig}}(\lambda)
\times \\
&{}_4F_3 (-n, n+\sig_0 +\delta_0 +1, \frac{\sig -\rho +\lam }2,
\frac{\sig -\rho -\lam }2; \sig_0 +1, \frac \sig 2,
\frac{\sig+1-\iota}2; 1),
\end{split}
\end{equation*}
\begin{equation*}
\begin{split}
&\qquad \mathcal F_{\pm1}(Q^{(0)}_{n, -1})
=-\frac{(\sig_0 +1)_n}{n!\sig}
\widetilde{w_{\sig}}(\lambda) (\pm\lam+\rho)
\times \\
&{}_4F_3 (-n ,  n+\sig_0 +\delta_1 +1, \frac{\sig -\rho +\lam }2,
\frac{\sig -\rho -\lam }2; \sig_0 +1, \frac{\sig} 2 +1,
\frac{\sig+1-\iota}2; 1).
\end{split}
\end{equation*}
\end{theo+}
\begin{proof}
The first result follows immediately
by using the Rodrigue's type formula.
 Indeed as in the proof of Theorem \ref{Rodrigue}
the function $Q_{n, 1}$ is a linear
combination of the functions
$w_{\sig+ 2m}(t)$,  whose
 spherical transform is, by
(\ref{eq:bern-var2}),
$$
\widetilde {w_{\sig+ 2m}}(\lam)=\frac{B_{m, \sig}(\lam)\widetilde{w_{\sig}}(\lambda)}{b_{m, \sig}}.
$$
Here we use the equation (2.3) and the formal self-adjointness
of $B_{m, \sig}(D)$ acting on functions of the form $\cosh^{-\sig-2j}t$,
which can be easily justified by direct computation or by the Plancherel
formula,
 since the function
$\widetilde{w_{\sig}}(i\lam)$ has exponential decay for $\lam \to \infty$.
 Our first formula follows then
 immediately by rewriting
the sum as a hypergeometric function.
For the second result we need to compute  the Cherednik-Opdam transform
of $(D+\rho)B_{m, \sig}(D) w_{\sig}(t)$. By the inversion
formula \cite[Theorem 9.13]{Opdam-acta}, we have,
$$
B_{m, \sig}(D) w_{\sig}(t)
=\int_{i\mathbb R^+}\big(\widetilde{w_{\sig}}(\lam)
G(\lam, t) +
\widetilde{w_{\sig}}(\lam)
G(-\lam, t)\big)
d\hat \mu(\lam).
$$
We let $D+\rho
$ acts on both
side,
$$
(D+\rho)B_{m, \sig}(D) w_{\sig}(t)
=\int_{i\mathbb R^+}\big(\widetilde{w_{\sig}}(\lam)
(D+\rho)
G(\lam, t) +
\widetilde{w_{\sig}}(\lam)
(D+\rho)G(-\lam, t)\big)
d\hat \mu(\lam).
$$
Rewriting the integral using (\ref{eig-eq})
we have
$$
(D+\rho)B_{m, \sig}(D) w_{\sig}(t)\int_{i\mathbb R^+}
\big(\widetilde{w_{\sig}}(\lam)
(\lam+\rho)
G(\lam, t) +
\widetilde{w_{\sig}}(\lam)
(-\lam+\rho)
G(-\lam, t)\big)d\hat \mu(\lam),
$$
which is equivalent to
that
$$
\mathcal F_{\pm1}\big((D+\rho)B_{m, \sig}(D) w_{\sig}\big)
(\lam)
=\widetilde{w_{\sig}}(\lam)
(\pm \lam+\rho).
$$
The rest follows by elementary computations.

\end{proof}
The Plancherel formula
for the transform $\mathcal F$
implies that
\begin{coro+} The functions $\{\mathcal F Q_{n,\pm 1}^{(0)}\}$
form an orthogonal basis of
$L^2(i\mathbb R^+, d\hat (\lam))\otimes \mathbb C^2$, and
their norms are given by
$$
\Vert \mathcal F Q_{n,\pm 1}^{(0)}\Vert =\Vert Q_{n,\pm
1}^{(0)}\Vert.$$
\end{coro+}

The general case of $k\ge 0$ can then
be obtained
by the same method using Proposition \ref{bern-k}.

\begin{theo+} Suppose $\sig > 2(\iota +b) $.  The Cherednik-Opdam transforms of the
Jacobi type functions $Q^{(k)}_{n, 1}$ and $Q^{(k)}_{n, -1}$ are given
by
\begin{equation*}
\begin{split}
&\quad \mathcal F_{\pm 1}(Q^{(k)}_{n, 1})
=\frac{(\sig_0+1)_n)}{n!} \widetilde{w_{\sig}}(\lambda)
\times \\
& \qquad
\sum_{m=0}^n\frac{(-n)_m (n+\sig_0 +\delta_0 +2k +1)_m}
{(\sig_0+1)_m m!} 
\frac{(\frac{\sig -\rho +\lam }2)_m (\frac{\sig -\rho -\lam
}2)_m} {(\frac{\sig}2)_m (\frac{\sig +1-\iota}2)_m} L_{k, m}(\lam)
\end{split}
\end{equation*}
and
\begin{equation*}
\begin{split}
&\quad \mathcal F_{\pm 1}(Q^{(k)}_{n, -1})= -\frac{(\sig_0
+1)_n}{n! \sig}
(\pm \lam+\rho) \widetilde{w_{\sig}}(\lambda)\times \\
&\qquad \sum_{m=0}^n\frac{(-n)_m (n+\sig_0 +\delta_0 +2k +1)_m}
{(\sig_0+1)_m m!} 
 \frac{(\frac{\sig -\rho +\lam }2)_m (\frac{\sig -\rho -\lam
}2)_m} {(\frac{\sig}2)_m (\frac{\sig +1-\iota}2)_m} M_{k, m}(\lam),
\end{split}
\end{equation*}
where $L_{k, m}$ and  $M_{k, m}$ are
the polynomials in (\ref{L-km})
and (\ref{M-km}).
For each fixed $k\ge 0$ the $\mathbb C^2$-valued
polynomials $\mathcal F(Q^{(k)}_{n, -1})(\lam)(\mathcal F_1(Q^{(k)}_{n, -1}(\lam), \mathcal F_1(Q^{(k)}_{n, -1}(\lam))$
 form  an orthogonal basis for
the space $L^2(i\mathbb R^+, d\hat \mu)\otimes \mathbb C^2$
and 
$$
\Vert \mathcal F(Q^{(k)}_{n, -1})\Vert \Vert Q^{(k)}_{n, -1})\Vert.
$$
\end{theo+}

\begin{rema+}
The functions $\mathcal F_{\pm 1}(Q^{(k)}_{n, \pm 1})$
are, apart from the common factor $\widetilde {w_\sig}(\lam)$,
 polynomials
of degree $2n+2k$ or $2n+2k+1$. However
 they are not the orthogonal polynomials
obtained by the usual Gram-Schmidt procedure
by ordering the monomials $1, \lam, \lam^2, \cdots$ according to
the degrees. Thus we have discovered some new
orthogonal polynomials of hypergeometric type.
\end{rema+}

\begin{rema+}
As an application, we recall that  \cite{peng-xu} 
we can find
an orthogonal basis
$\{P_n^{(\alpha,k)}(2|z|^2-1)|z|^k\}_{n=0}^\infty$,
for each fixed $k$, 
of the $L^2$-space of radial
functions
on the unit disk with the invariant
measure $d\iota(z)=(1-|z|^2)^{-2}dz\wedge d\bar z$,
as a symmetric space $SU(1, 1)/U(1)$,
They are just the functions $Q^{(k)}_{n, \pm 1}(t)$ in variable
$\tanh t=|z|$. So Theorem 4.3 gives their
Helgason transform \cite{Helg-2} on the unit disk.
\end{rema+}

Finally we mention that the result obtained in
this paper can also be used to find a family
of  polynomials of Wilson type  
orthogonal with respect to the so-called
asymmetric Harish-Chandra  $c$-function
\cite{Cherednik-imrn97}, see also 
\cite{WG-05}.

\newcommand{\noopsort}[1]{} \newcommand{\printfirst}[2]{#1}
  \newcommand{\singleletter}[1]{#1} \newcommand{\switchargs}[2]{#2#1}
  \def\cprime{$'$} \def\cprime{$'$}
\providecommand{\bysame}{\leavevmode\hbox
to3em{\hrulefill}\thinspace}
\providecommand{\MR}{\relax\ifhmode\unskip\space\fi MR }
\providecommand{\MRhref}[2]{%
  \href{http://www.ams.org/mathscinet-getitem?mr=#1}{#2}
} \providecommand{\href}[2]{#2}

\end{document}